\documentclass[number,citesort,seceqn,MSNbibl,dvips]{arxbj}
\usepackage{upgreek}

\aid{0}
\volume{17}
\issue{2}
\pubyear{2011}
\firstpage{609}
\lastpage{627}
\doi{10.3150/10-BEJ294}

\makeatletter
\newtheorem{theorem}{Theorem}[section]
\newtheorem{lemma}{Lemma}[section]
\newtheorem{coro}{Corollary}[section]

\newremark{remark}{Remark}[section]
\newremark{example}{Example}[section]

\def\sklfrac#1#2{#1/#2}

\def\afrac#1#2{#1/(#2)}

\newcommand{\implies}{\Longrightarrow}
\newcommand{\varlimsup}{\overline{\mathop{\mathrm{lim}}}}
\newcommand{\Real}{\mathbb{R}}
\newcommand{\Z}{\mathbb{Z}}
\newcommand{\eps}{\varepsilon}
\newcommand{\s}{\mathcal{S}}
\newcommand{\one}[1]{\mathbf{1}_{\{#1\}}}
\renewcommand{\P}{\mathbb{P}}
\newcommand{\E}{\mathbb{E}}
\newcommand{\argmin}{\mathop{\mathrm{argmin}}}
\newcommand{\eqref}[1]{(\ref{#1})}
\newcommand{\argmax}{\mathop{\mathrm{argmax}}}
\newcommand{\Mu}{\mathcal{M}}
\newcommand{\Partial}{\mathcal{D}_{12}}

\makeatother

\begin{document}
\begin{frontmatter}

\title{On the Viterbi process with continuous\break\ state space}
\runtitle{On the Viterbi process with continuous state space}

\begin{aug}
\author[1]{\fnms{Pavel} \snm{Chigansky}\corref{}\ead[label=e1]{pchiga@mscc.huji.ac.il}}
\and
\author[1]{\fnms{Yaacov} \snm{Ritov}\ead[label=e2]{yaacov.ritov@gmail.com}}

\runauthor{P. Chigansky and Y. Ritov}

\address[1]{Department of Statistics,
The Hebrew University,
Mount Scopus, Jerusalem 91905,
Israel.\\
\printead{e1,e2}}
\end{aug}

% HISTORY:
\received{\smonth{9} \syear{2009}}
\revised{\smonth{3} \syear{2010}}

% ABSTRACT
%
\begin{abstract}
This paper deals with convergence of the maximum a posterior probability
path estimator in hidden Markov models. We show that when the state
space of the hidden process is continuous,
the optimal path may stabilize in a way which is essentially different
from the previously considered finite-state setting.
\end{abstract}

% KEYWORDS
%

\end{frontmatter}

%s1 ###
\section{Introduction}
Consider a standard hidden Markov model $(X,Y)$, where
$X=(X_n)_{n\in\Z_+}$ and $Y=(Y_n)_{n\in\Z_+}$ are the \textit{hidden}
state and the \textit{observation} processes, respectively. The state
process $X$ is Markov with values in a subset $\s\subseteq\Real$,
transition probability
$Q$ and initial distribution $\Mu$: for all measurable subsets
$A\subseteq\s$,
\begin{eqnarray*}
 \P(X_1\in A)&=&\Mu(A),\\
 \P(X_n \in A |X_{n-1})&=&Q(X_{n-1},A), \qquad \P\mbox{-a.s., }  n>1.
\end{eqnarray*}
We shall consider either countable $\s$, in which case $q(u,v):=Q(u,\{
v\})$ and $\mu(u):=\Mu(\{u\})$, or
$\s=\Real$, assuming that $Q(u,\mathrm{d}v)$ and $\Mu(\mathrm{d}u)$ have densities
$q(u,v)$ and $\mu(u)$ with respect to the Lebesgue measure. The precise meaning
of $q(u,v)$ and $\mu(u)$ should be obvious from the context.

The observed process $Y$ forms a sequence of conditionally independent
random variables,
given $X_{1:\infty}=(X_1,X_2,\ldots)$, with the \textit{observation}
density $p$:
\[
\P(Y_n\in B|X_{1:\infty}) = \int_B p(X_n,y)\,\mathrm{d}y,\qquad \P\mbox{-a.s.},
\]
for any Borel $B\subseteq\Real$.

The path estimation problem is to reconstruct the trajectory of the
hidden process\footnote{Hereafter, for $x\in\Real^n$, $x_m$ stands for
the $m$th entry of $x$ and $x_{k:m}$, $k\le m,$
denotes the vector $x=(x_k,\ldots,x_m)$; $|x_{1\dvtx n}|=\max_i |x_i|$ and
$\|x_{1\dvtx n}\|=\sqrt{\sum_{i=1}^n x_i^2}$.
} $X_{1\dvtx n}=(X_1,\ldots,X_n)$,
given the realization of $Y_{1\dvtx n}=(Y_1,\ldots,Y_n)$ for a fixed horizon
$n\ge1$. If $\s$ is a discrete set, a natural estimator is the
maximizer of the a posterior probability (MAP estimator):
\[
\hat X^n_{1\dvtx n}:=\argmax_{x_{1\dvtx n}\in\s^n}\P
(X_{1\dvtx n}=x_{1\dvtx n}|Y_{1\dvtx n} ),
\]
where the optimal path is chosen according to the
lexicographical order on $\s^n$, induced by an order on $\s$, whenever
the maximum is not unique.
The obtained path minimizes the probability of error among all
estimators depending on $Y_{1\dvtx n}$, that is,
\[
\P(\hat X^n_{1\dvtx n}\ne X_{1\dvtx n})\le\P(\xi_{1\dvtx n}\ne X_{1\dvtx n})
\qquad \mbox{for all $\sigma\{Y_1,\ldots,Y_n\}$-measurable $\xi_{1\dvtx n}$}.
\]

By Bayes' formula,
\[
\P (X_{1\dvtx n}=x_{1\dvtx n}|Y_{1\dvtx n} ) = \frac{L_n(x_{1\dvtx n},Y_{1\dvtx n})}{\sum
_{u_{1\dvtx n}\in\s^n}L_n(u_{1\dvtx n},Y_{1\dvtx n})},
\]
where $L_n$ is the ``posterior'' likelihood:
%
%e1 ###
\begin{equation}
\label{Ln}L_n(x_{1\dvtx n};y_{1\dvtx n}) = \mu(x_1)p(x_1,y_1)\prod_{m=2}^n
q(x_{m-1},x_{m})p(x_m,y_m), \qquad  x_{1\dvtx n}\in\s^n,
\end{equation}
and hence
\[
\hat X^n_{1\dvtx n}=\argmax_{x_{1\dvtx n}\in\s^n}L_n(x_{1\dvtx n},Y_{1\dvtx n}).
\]

Due to the product structure of $L_n$, the search for the maximizing
path can be carried out efficiently by a
dynamic programming procedure, called the \textit{Viterbi algorithm},
after A.~Viterbi, who introduced it in the context of
error correction codes.

When the next\vspace*{1pt} observation, $Y_{n+1}$, is added, the optimal path may
change entirely, that is, for any $m=1,\ldots,n$,
$\hat X^{n+1}_{1:m}$ is, in general, different from $\hat X^{n}_{1:m}$.
In practical terms, the latter
means that\footnote{$\# A$ stands for the cardinality of a set $A.$} $\#
\s$ optimal path candidates of length $n$ are to be kept in memory at
each time $n$. This
motivates the question of whether the optimal path stabilizes as the
number of observations grows to infinity or, more precisely,
whether the limit
%
%e2 ###
\begin{equation}
\label{hat}
\hat X_{1:m}=\lim_{n\to\infty}\hat X^{n}_{1:m}
\end{equation}
exists $\P$-a.s.~for each fixed $m\ge1$. If such a limit exists, it
defines a random process with paths in $\s^\infty$, named (in
\cite{LeKo09}) \textit{the Viterbi process}.

An affirmative answer to this question was given in \cite{CaRo02} (see
also \cite{Ko96}) under a sufficient condition (see \eqref{CR} below)
which also ensures that the limit sequence $\hat X=(\hat X_m)_{m\ge1}$ is
a regenerative process. More precisely, a sequence of stopping times
can be constructed (see \cite{Ca06}), splitting the process
$\hat X$ into cycles that are i.i.d.~and independent of the initial
delay. In particular, by the regenerative property, $\hat X$ satisfies
the classical limit laws, such as the law of large numbers (LLN) and
the central limit theorem (CLT).

In fact, the existence of such renewal times under the condition \eqref
{CR} can be deduced by a simple argument
(reproduced, for completeness, in Section~\ref{sec-2}). A more delicate
construction in \cite{LeKo08,LeKo09}
verifies \eqref{hat} under conditions weaker than \eqref{CR}.

In this paper, we revisit the question of the existence of the limit
\eqref{hat} for hidden Markov models (HMMs) with continuous state
spaces, that is, when $\s=\Real$ and
for each $u\in\Real$, the transition kernel $Q(u,\mathrm{d}v)$ and the initial
distribution $\Mu(\mathrm{d}v)$ have densities $q(u,v)$ and $\mu(v)$, respectively,
with respect to the Lebesgue measure. By Bayes' formula,
the conditional law of the vector $X_{1\dvtx n}$ given $Y_{1\dvtx n}$ has the
density $\psi_n$ with respect to the Lebesgue measure on $\Real^n$:
\[
\psi_n(x_{1\dvtx n}):=\frac{L_n(x_{1\dvtx n};Y_{1\dvtx n})}{\int_{\Real
^n}L_n(u_{1\dvtx n};Y_{1\dvtx n})\,\mathrm{d}u_1 \cdots \,\mathrm{d}u_n},
\]
with $L_n$ defined as in \eqref{Ln}.
The MAP path estimator is
\[
\hat X^n_{1\dvtx n} := \argmax_{x_{1\dvtx n}\in\Real^n}\psi_n(x_{1\dvtx n})=\argmax
_{x_{1\dvtx n}\in\Real^n}L_n(x_{1\dvtx n};Y_{1\dvtx n}),
\]
where, as in \eqref{hat}, the maximum is chosen according to the
lexicographical order on
$\Real^n$ (induced, e.g., by $<$ on $\Real$) in case of ambiguity.

Note that for any $\sigma\{Y_1,\ldots,Y_n\}$-measurable random vector
$\xi_{1\dvtx n}$ and $\eps>0$,
\[
\P(|X_{1\dvtx n}- \xi_{1\dvtx n}|\le \eps)=
\E\P(|X_{1\dvtx n}- \xi_{1\dvtx n}|\le \eps|Y_{1\dvtx n} )=
\E\int_{[-\eps,\eps]^n} \psi_n(x_{1\dvtx n}+\xi_{1\dvtx n})\,\mathrm{d}x_1\cdots \,\mathrm{d}x_n
\]
and hence the estimator $\hat X^n_{1\dvtx n}$ is optimal in the sense that
\begin{eqnarray*}
\lim_{\eps\to0}\eps^{-n}\P(|X_{1\dvtx n}- \xi_{1\dvtx n}|\le \eps)&=&
\E\psi_n(\xi_{1\dvtx n})\le\E\max_{x_{1\dvtx n}\in\Real^n}\psi_n(x_{1\dvtx n})\\
&&{}=\lim_{\eps\to0}\eps^{-n}\P( |X_{1\dvtx n}- \hat X^n_{1\dvtx n} |\le\eps)
\end{eqnarray*}
whenever interchanging the expectation and the limit is possible.
Roughly, this means that $\hat X^n_{1\dvtx n}$ yields the best ``small''
credible intervals among all other path
estimates.\footnote{In fact, this optimality interpretation turns out to
be meaningful even in the infinite-dimensional function space; see \cite
{ZD87,ZD88}.}

As in state estimation problems such as filtering, the exact
calculation of $\hat X^n_{1\dvtx n}$ is impossible beyond a number
of models with a special structure, most notably Kalman's linear
Gaussian setting. A number of efficient numerical techniques, such
as particle filters, have been developed (see, e.g., \cite{CaMoRy05})
to approximate the conditional law of the hidden state process.
In this paper, we are concerned with the convergence properties of the
MAP paths, leaving the computational
issues for further investigation.

In Section \ref{sec-2}, we explore, through a number of examples,
various patterns of convergence encountered in \eqref{hat}, when
the hidden state space is continuous. We also give an example of HMM,
for which the MAP path does not converge as the estimation time
horizon increases.
In Section \ref{sec-3}, we prove a more general result, deducing the
existence of the limit \eqref{hat} from certain strong
$\log$-concavity of the transition and observation densities. The
\hyperref[app-A]{Appendix} contains a lemma which is used in the proof of the main
result and
may be of independent interest. Finally, a short discussion of the
results appears in Section \ref{sec-4}.

%s2 ###
\section{Examples}\label{sec-2}

Let us briefly recall the essential elements of the proof in the finite setting $\s=\{1,\ldots,d\}$.
For simplicity, consider an irreducible finite (and thus recurrent)
chain $X$ and define
\[
D_i = \{y\in\Real\dvtx  q(x_1,i)p(i,y)q(i,x_3)>
q(x_1,x_2)p(x_2,y)q(x_2,x_3), \forall x_2\ne i, x_1, x_3\in\s\}.
\]
Suppose that, for a pair of states $j_0$ and $i_0$,
%
%e3 ###
\begin{equation}\label{CR}
\int_{D_{i_0}} p(j_0,y)\,\mathrm{d}y>0.
\end{equation}
Recall the definition of $L_n$ in \eqref{Ln} and note that on the event
$A_m=\{X_m =j_0, Y_m\in D_{i_0}\},$ with a fixed $m>1$ and all $n>m,$
\begin{eqnarray*}
L_n(x_{1\dvtx n},Y_{1\dvtx n})
&=&L_{m-1} (x_{1:m-1},Y_{1:m-1} )
\\
&&{}\times q(x_{m-1},x_m)p(x_m,Y_m)q(x_m,x_{m+1})L_{m+1,n}\bigl(x_{(m+1):n},Y_{(m+1):n}\bigr)
\\
&\le& L_{m-1} (x_{1:m-1},Y_{1:m-1} )
\\
&&{}\times q(x_{m-1},i_0)p(i_0,Y_m)q(i_0,x_{m+1}) L_{m+1,n}\bigl(x_{(m+1):n},Y_{(m+1):n}\bigr)
\end{eqnarray*}
for an appropriate function $L_{m+1:n}$ and where equality is attained
only at a path $x_{1:m}$ with $x_m=i_0$. Hence, the $m$th entry of the
optimal path must equal $i_0$
for any $n\ge m$, that is, $\hat X^n_{m}=i_0$.
But, then, given $\hat X_m^n$, the first $m$ entries of the optimal
path depend only on the values of $Y_1,\ldots,Y_m$ and are not
affected by $Y_k$, $k>m$. Hence, the limit \eqref{hat} exists on the
event $A_m$. Since the chain $(X,Y)$ is recurrent, for any
fixed $m$, one of the events $A_{m+1},A_{m+2},\ldots$ occurs $\P
$-a.s.~and thus \eqref{hat} holds $\P$-a.s.

Using the same basic idea, let $\tau(k)$, $k\ge0,$ be the times at
which the chain $(X,Y)$ revisits the set $\{j_0\}\times D_{i_0}$:
\begin{eqnarray*}
\tau(0)&=&1,\\
\tau(k) &=& \inf\{n> \tau(k-1)\dvtx  X_n = j_0, Y_n \in D_{i_0}\},\qquad k\ge1.
\end{eqnarray*}
By construction, for any $k$, on the event $\{\tau(k)\le n\}$,
\[
L(x_{1\dvtx n};Y_{1\dvtx n}) \le L(x^\tau_{1\dvtx n};Y_{1\dvtx n})\qquad \forall x_{1\dvtx n}\in\s^n,
\]
where $x^\tau_{1\dvtx n}$ is the vector which coincides with $x_{1\dvtx n}$ at
all but the indices
$\tau(1),\ldots,\tau(k)$, where its entries equal $i_0$.

The upper bound is attained if $L(x_{1\dvtx n};Y_{1\dvtx n})$ is maximized over
$x_{1\dvtx n}$, constrained to
$x_{\tau(1)}=\cdots=x_{\tau(k)}=i_0$. Since each $x_{\tau(\ell)}$, $\ell
=1,\ldots,k$, appears in the product $L(x_{1\dvtx n};Y_{1\dvtx n})$ in
three adjacent terms, the optimal choice for each segment $x_{\tau(\ell
-1)+1:\tau(\ell)-1}$, $\ell=1,\ldots,k$, is determined only by the
values of
$Y_{\tau(\ell-1)+1},\ldots,Y_{\tau(\ell)-1}$. Hence, in particular, the limit
$\lim_{n\to\infty}\hat X^n_{1:m}$ exists on any of the events $\{\tau
(k-1)<m\le\tau(k)<\infty\}$, $k\ge1$.
By recurrence of $j_0$ and the condition~\eqref{CR}, $\P(\tau(k)<\infty
)=1$ and
$\lim_{k\to\infty}\tau(k)=\infty$, $\P$-a.s., which verifies the
existence of the limit \eqref{hat}.

The stopping times $\tau(k)$, $k\ge1$, form a renewal process, with
respect to which both $(X,Y)$ and $\hat X =(\hat X_m)_{m\ge1}$
are regenerative (see \cite{Ca06} for more details). As pointed out in
\cite{LeKo08},
the condition~\eqref{CR} can be quite restrictive, especially when the
transition matrix is sparse.
The convergence in \eqref{hat} and the regenerative property are
verified in \cite{LeKo08} under less conservative conditions,
using a more sophisticated construction of the renewal times.

In summary, both \cite{CaRo02} and \cite{LeKo08} deduce the existence
of the limit in \eqref{hat} from the explicit construction of stopping times,
based on the discreteness of the hidden process state space. The
following example shows that this still may be possible in HMMs with
continuous state spaces.
%ex2.1
\begin{example}
Consider a linear HMM with Laplacian state and Gaussian observation noises:
\[
\mu(u)= \frac1 4 \mathrm{e}^{-|u|/2},\qquad
q(u,v)=\frac1 4 \mathrm{e}^{-|u-v|/2},\qquad
p(x,y) =\frac1 {\sqrt{2\uppi}}\mathrm{e}^{-(x-y)^2/2}.
\]
In this case, the MAP path is given by
\[
\hat X^n_{1\dvtx n} =\argmin_{x_{1\dvtx n}\in\Real^n} \Biggl(|x_1| + (x_1-Y_1)^2+ \sum
_{m=2}^n|x_{m-1}-x_m|+ (x_m-Y_m )^2 \Biggr).
\]
Consider the function $x\mapsto f(x):=|a-x|+(x-y)^2+|x-b|$ for fixed
$a,b,y\in\Real$.
Suppose, without loss of generality, that $a\le b$ and note that $f$,
being strictly convex, is minimized at a
unique point $x^*=\argmin_{x\in\Real}f(x)$. If $y\in[a,b]$, then,
clearly, $x^*\in[a,b]$ and since
$f(x)=-a+(y-x)^2+b$ on this interval, we have $x^*=y$. Consider the
case $y\le a$ and suppose $x^*< a$. For $x< a$,
$f(x)=a-x +(y-x)^2+b-x$ and hence $x^*=y+1$. By strict convexity, this
implies that $x^*=y+1$ if $y<a-1$ and
that $x^*\ge a$ otherwise. Clearly, $x^*\le b$, that is, $x^*\in
[a,b]$, which, in turn, implies that $x^*=a$ for $y\in[a-1,a)$.
Similar calculations reveal that $x^*=y-1$ if $y>b+1$ and $x^*=b$ if
$y\in(b,b+1]$.

To summarize, $x^*\in[y-1,y+1]$ for any $a,b,y\in\Real$ and $x^*=y$,
whenever $a\le y\le b$. In particular,
$\hat X^n_{m-1}\in[Y_{m-1}-1,Y_{m-1}+1]$ and $\hat X^n_{m+1}\in
[Y_{m+1}-1,Y_{m+1}+1]$
for any $n\ge m+1$.
Hence, on the event
\[
A_m:=\{Y_{m-1}+1\le Y_m\le Y_{m+1}-1\},
\]
$Y_m\in[\hat X^n_{m-1}, \hat X^n_{m+1}]$ and, consequently, $\hat X^n_m=Y_m$.
This, in turn, implies that $\hat X^n_{1:m}=\hat X^{m+1}_{1:m}$ for all
$n\ge m+1$
and the existence of the limit \eqref{hat} on any of $A_k$, $k\ge m+1$.
Clearly, the $A_k$'s occur infinitely often and hence,
as in the discrete case, $\hat X^n_{1:m}$ ceases to change, starting
from some random, but $\P$-a.s. finite, time $n$. In particular,
\eqref{hat} holds $\P$-a.s.
\end{example}

However, splitting the optimal trajectory into unrelated segments is
not the only way to get the convergence in \eqref{hat}: the following
example shows that
the limit may exist without ever being actually attained.
%ex2.2
\begin{example}
Consider the linear Gaussian HMM with
\[
\mu(u)=\frac1{\sqrt{2\uppi}} \mathrm{e}^{-u^2/2},\qquad
q(u,v) = \frac1{\sqrt{2\uppi}}
\mathrm{e}^{-(u-v)^2/2},\qquad
p(x,y)=\frac1{\sqrt{2\uppi}} \mathrm{e}^{-(x-y)^2/2}.
\]
In this case, the conditional law of $X_{1\dvtx n}$, given $Y_{1\dvtx n}$, is
Gaussian and hence
\[
\hat X^n_{1\dvtx n}=\E(X_{1\dvtx n}|Y_{1\dvtx n} ).
\]
For any fixed $m\ge1$, the process $\hat X^n_{1:m}=\E
(X_{1:m}|Y_{1\dvtx n})$, $n\ge m$, is a uniformly integrable vector-valued
martingale and hence
the limit \eqref{hat} exists by the martingale convergence. In fact,
Kalman linear filtering theory (see, e.g., \cite{KS72})
guarantees that in this case (of controllable and observable dynamics)
the stronger $\P$-a.s.~exponential convergence holds (see also Remark
\ref{rem-G} below).

Moreover, $\E(X_{1:m}|Y_{1\dvtx n})$ is a deterministic linear map of
$Y_{1\dvtx n}$ and a
calculation reveals that it actually depends on each one of the
components in $Y_{1\dvtx n}$. Since $Y_{1\dvtx n}$ is a non-degenerate Gaussian vector,
\[
\P(\hat X^n_j=\hat X^{n'}_j, \mbox{for some\ } j \le m )=0
\]
for any $n'>n\ge m$.
\end{example}

Finally, the next example demonstrates that a finite limit in \eqref
{hat} may not exist, even when the hidden state chain is
positive recurrent and has countably many states. In fact, it also
shows that the optimal MAP path may not be an adequate estimate:
in this case, a trajectory of a positive recurrent chain $V$ is
estimated as a constant trajectory, diverging to infinity, as $n\to
\infty$.
%ex2.3
\begin{example}
Consider the HMM with the hidden state process $X_n = (U_n,V_n)$,
consisting of independent
components $U$ and $V$. The process $U=(U_n)_{n\ge1}$ is a sequence of
i.i.d.~random variables uniformly distributed over $[0,1]$.

$V=(V_n)_{n\ge1}$ is a random walk on positive integers with
reflecting boundary at $\{1\}$ and the
transition probabilities $P(1,1)=1-\eps$, $P(1,2)=\eps$ and, for $i\ge2$,
%
%e4 ###
\begin{equation}\label{Pij}
P(i,j) =
\cases{
\displaystyle \eps\frac{ (\afrac{i}{i+1} )^2}{1+ (\afrac{i}{i+1} )^2}, &\quad $j = i+1$,\cr
\displaystyle \eps\frac{1}{1+ (\afrac{i}{i+1} )^2}, &\quad $j =i-1$,\cr
1-\eps, &\quad $j=i$,
}
\end{equation}
where $\eps>0$ is a small fixed constant (in fact, we shall later
choose $\eps<\mathrm{e}^{-2}/(1+\mathrm{e}^{-2})=0.119\ldots$). $V$ is a positive
recurrent Markov chain with the unique invariant distribution
%
%e5 ###
\begin{equation}\label{piexp}
\pi(j)=
\cases{
 \frac1 5 C \bigl(1+ \bigl(\frac{1}{2} \bigr)^2 \bigr), &\quad$j=1$, \cr
\displaystyle C\frac1 {j^2} \biggl(1+ \biggl(\frac{j}{j+1} \biggr)^2 \biggr), &\quad $j> 1$,
}
\end{equation}
where $C$ is the normalization constant, independent of $\eps$. We
shall assume that $V$ is stationary,
that is, it is started from $V_1\sim\pi$. Stationarity is not really
required in what follows and is solely a matter
of aesthetics (e.g., $\P(V_1=j)= C/j^2$ will work as well).\vspace*{1pt}

Let $a_0=0$, $a_i=8\sum_{j=1}^i(1/9)^j$, $i=1,2,\ldots,$ and set $A_i =
[a_{i-1},a_i)$, $i\ge1$. Denote by $\ell_i=8(1/9)^i$ the length of the interval
$A_i$ and note that $[0,1)=\bigcup_{i=1}^\infty A_i$.

Now, consider the observation density
\[
p ((u,v),y )= \one{y\in[0,1]}\one{u\notin\bigcup_{i=1}^v A_i}+\sum
_{i=1}^v \ell_i^{-1}\one{(u,y)\in A_i\times A_i}.
\]

As we show below, the MAP estimates of $U_{1\dvtx n}$ and $V_{1\dvtx n}$ are
given by\footnotemark:
%
%e6 ###
\begin{eqnarray}
\label{opt}
\hat U^n_m &=&
\sum_{j=1}^\infty a_{j-1} \one{Y_m\in A_j},\qquad
m=1,\ldots,n,\nonumber
\\[-8pt]\\[-8pt]
\hat V^n_m &=&
\cases{
2, &\quad $j^*(n)=1$,\cr
j^*(n), &\quad $j^*(n)>1$,
}\nonumber
\end{eqnarray}
where
$j^*(n):=\max\{j\dvtx \sum_{k=1}^n \one{Y_k\in A_j}>0 \}$.
\footnotetext
{
The choice of $\hat U^n_m$ is not unique, unless the lexicographic
order is imposed:
for example, $\hat U^n_m:=Y_m$ yields the same value of the likelihood.
}
Since all $A_j$'s have positive Lebesgue measure, $j^*(n)\nearrow\infty
$ as $n\to\infty$ and, consequently,
for any fixed $m\ge1$,
\[
\lim_{n\to\infty}\hat V^n_m =\lim_{n\to\infty}j^*(n)=\infty,\qquad \P\mbox{-a.s.}
\]

Before proving \eqref{opt}, we shall briefly explain why the optimal
path of such a form should be
anticipated. Note that since $U_i$'s are uniformly distributed in
$[0,1]$, the choice
of $\hat U^n_i$'s influences the likelihood \eqref{Ln} only through the
observation densities. More precisely,
whenever $\{Y_m\in A_i\}$ is observed, the maximal gain of $\ell
^{-1}_i$ is obtained if $\hat U^n_m\in A_i$
and $\hat V^n_m\ge i$ are chosen. On the other hand, the transition
probabilities of \eqref{Pij} favor
paths $\hat V^n_{1\dvtx n}$ without jumps. Hence, the optimal path $\hat
V^n_{1\dvtx n}$ should be constant and
large enough to allow access to the narrowest $A_i$ visited by $Y_m$'s
so far, that is, greater or equal to
$j^*(n)$. However, if constant $\hat V^n_{1\dvtx n}$ is chosen, it cannot be
too large, as this would decrease
the likelihood through the term $\pi(\hat V^n_1)$, due to the fast tail
decay of the initial distribution $\pi$.
This heuristics is implemented by an appropriate balancing between all
the ingredients of the model.

We shall first check \eqref{opt} in the case $j^*(n)>1$. To this end,
consider the ratio
%
%e7 ###
\begin{equation}
\label{show0}
\frac
{
L_n ((u_{1\dvtx n},v_{1\dvtx n}),Y_{1\dvtx n} )
}
{
L_n ((\hat U^n_{1\dvtx n},\hat V^n_{1\dvtx n}),Y_{1\dvtx n} )
}=
\frac
{
\pi(v_1)
}
{
\pi(j^*(n))
}
\prod_{m=2}^n \frac
{
P (v_{m-1},v_m )
}
{
P (j^*(n),j^*(n) )
}
\prod_{m=1}^n
\frac
{
p ((u_m,v_m),Y_m )
}
{
p ((\hat U^n_m,j^*(n)),Y_m )
}
\end{equation}
for an arbitrary $u_{1\dvtx n}$ and $v_{1\dvtx n}$.
Let $N$ be the number of jumps in $v_{1\dvtx n}$ and $v^*(n)=\max_{k=1,\ldots,n}v_k$.
Note that $P (v_{m-1},v_m )=1-\eps$ when $v_{m-1}= v_m$ and
$P (v_{m-1},v_m )\le\eps$ otherwise. Hence, as $P (j^*(n),j^*(n)
)=1-\eps$,
\[
\prod_{m=2}^n \frac
{
P (v_{m-1},v_m )
}
{
P (j^*(n),j^*(n) )
}\le \biggl(\frac{\eps}{1-\eps} \biggr)^N.
\]
Further, note that on the event $\{Y_m\in A_i\}$,
$
p ((u_m,v_m),Y_m ) \le1 \vee\ell^{-1}_{i}=\ell^{-1}_{i}
$
and
$p ((\hat U^n_m,j^*(n)), Y_m )=\ell^{-1}_{i}$, thus
\[
\frac
{
p ((u_m,v_m),Y_m )
}
{
p ((\hat U^n_m,j^*(n)),Y_m )
}\le1.
\]
Moreover, on $\{Y_m\in A_{j^*(n)}\}$,
\[
\frac
{
p ((u_m,v_m),Y_m )
}
{
p ((\hat U^n_m,j^*(n)),Y_m )
}\le
\frac
{
\one{v^*(n)<j^*(n)}+\ell^{-1}_{j^*(n)}\one{v^*(n)\ge j^*(n)}
}
{
\ell^{-1}_{j^*(n)}
}
\le
\frac
{
\ell^{-1}_{v^*(n)\wedge j^*(n)}
}
{
\ell^{-1}_{j^*(n)}
}.
\]
Plugging these inequalities into \eqref{show0}, we get
%
%e8 ###
\begin{eqnarray}\label{show}
\frac
{
L_n ((u_{1\dvtx n},v_{1\dvtx n}),Y_{1\dvtx n} )
}
{
L_n ((\hat U^n_{1\dvtx n},\hat V^n_{1\dvtx n}),Y_{1\dvtx n} )
}
&\le&
\frac
{
\pi(v_1)
}
{
\pi(j^*(n))
}
\biggl(\frac
{
\eps
}
{
1-\eps
}
\biggr)^N
\frac
{
\ell^{-1}_{v^*(n)\wedge j^*(n)}
}
{
\ell^{-1}_{j^*(n)}
}\nonumber
\\[-8pt]\\[-8pt]
\nonumber&=&
\frac
{
\pi(v_1)
}
{
\pi(v^*(n))
}
\tilde\eps^N
\frac
{
\pi(v^*(n))
}
{
\pi(j^*(n))
}
\frac
{
\ell^{-1}_{v^*(n)\wedge j^*(n)}
}
{
\ell^{-1}_{j^*(n)}
},
\end{eqnarray}
where we define $\tilde\eps:=\eps/(1-\eps)$ for the purposes of brevity.
Since $N\ge v^*(n)-v_1$,
\begin{eqnarray*}
\frac
{
\pi(v_1)
}
{
\pi(v^*(n))
}
\tilde\eps^{N}
&\le&
\frac
{
\pi(v_1)
}
{
\pi(v^*(n))
}
\tilde\eps^{v^*(n)-v_1}
\le
\biggl(\frac{v^*(n)}{v_1} \biggr)^2
\frac
{
1+ (\afrac{v_1}{v_1+1} )^2
}
{
1+ (\afrac{v^*(n)}{v^*(n)+1} )^2
} \tilde\eps^{v^*(n)-v_1}
\\
&\le&
\biggl(\frac{v^*(n)}{v_1} \biggr)^2
\tilde\eps^{v^*(n)-v_1},
\end{eqnarray*}
where, in the second inequality, we have used the expression for $\pi
(j)$, $j>1,$ from \eqref{piexp}.
In fact, the inequality is also true for $v^*(n)=v_1=1,$ as both the
right- and left-hand sides become $1$, and
for $v^*(n)>v_1=1$, as $\pi(1)$ is less than $C\frac1 {j^2} (1+ (\frac
{j}{j+1} )^2 )$ evaluated at $j:=1$.

The function $x\mapsto x^2\tilde\eps^x$ attains its maximum at
$x^*=2/\log\tilde\eps^{-1}$ and is strictly decreasing on $(x^*,\infty)$.
Hence, with $\tilde\eps<\mathrm{e}^{-2}$, that is, with $\eps<
\mathrm{e}^{-2}/(1+\mathrm{e}^{-2})$, for any $y> x\ge1$, $(y/x)^2\tilde\eps^{y-x}<1$
and hence
%
%e9 ###
\begin{equation}
\label{raz}
\frac
{
\pi(v_1)
}
{
\pi(v^*(n))
}
\tilde\eps^{N}\le1.
\end{equation}
The equality holds if and only if $v_{1\dvtx n}$ is a constant path, that
is, $v_{m}=v^*(n)$ for all $m=1,\ldots,n$.

Further, if $v^*(n)\le j^*(n)$, then\vspace*{-2pt}
%
%e10 ###
%
\begin{eqnarray}\label{dva}
\frac
{
\pi(v^*(n))
}
{
\pi(j^*(n))
}
\frac
{
\ell^{-1}_{v^*(n)\wedge j^*(n)}
}
{
\ell^{-1}_{j^*(n)}
}&\le&
\biggl(\frac{j^*(n)}{v^*(n)} \biggr)^2
\frac
{
1+ (\afrac{v^*(n)}{v^*(n)+1} )^2
}
{
1+ (\afrac{j^*(n)}{j^*(n)+1} )^2
}
(1/9)^{j^*(n)-v^*(n)}\\
&\le&
\biggl(\frac{j^*(n)}{v^*(n)} \biggr)^2
(1/9)^{j^*(n)-v^*(n)}\le1,
\end{eqnarray}
where the latter inequality holds since $1/9< \mathrm{e}^{-2}/(1+\mathrm{e}^{-2})$.

The sequence $\pi(j)$ attains its unique maximum at $j:=2$ and is
strictly decreasing for $j\ge2$.
Hence, if $v^*(n)>j^*(n)\ge2$,
then\vspace*{-2pt}
\[
\frac
{
\pi(v^*(n))
}
{
\pi(j^*(n))
}
\frac
{
\ell^{-1}_{v^*(n)\wedge j^*(n)}
}
{
\ell^{-1}_{j^*(n)}
}< 1.
\]

Plugging \eqref{raz} and \eqref{dva} into \eqref{show} yields the
following inequality for any $u_{1\dvtx n}$ and $v_{1\dvtx n}$:\vspace*{-2pt}
\[
L_n ((u_{1\dvtx n},v_{1\dvtx n}),Y_{1\dvtx n} )
\le
L_n ((\hat U^n_{1\dvtx n},\hat V^n_{1\dvtx n}),Y_{1\dvtx n} ),
\]
which saturates if and only if $v_m=j^*(n)$, $m=1,\ldots,n$, thus
verifying the optimality of \eqref{opt}
on the event $\{j^*(n)>1\}$.

We shall omit the details in the case $\{j^*(n)=1\}$, which is treated
similarly: the optimal value $\hat V^n_m=2$
is obtained since $\pi(j)$ is maximal at $j=2$. Of course, as $j^*(n)$
eventually leaves the state~$1$,
the exact value is irrelevant for the main point of the present
example, that is, the divergence $\lim_{n\to\infty} \hat V^n_m=\infty$.
\end{example}

%s3 ###
\section{Convergence in the case of log-concave densities}\label{sec-3}
In this section, we establish the existence of the limit \eqref{hat},
deducing it from
certain strong $\log$-concavity properties of the densities $q$ and
$p$. Hereafter, the following assumptions are
in force:
\begin{enumerate}[(a1)]
\item[(a1)]\label{a0} the initial state density $\mu$ is a $C^2(\Real)$ $\log
$-concave function on $\Real$ and $-\log\mu(u)\ge0$;

\item[(a2)]\label{a1} the hidden state transition density $q$ is a $C^2(\Real
^2)$ $\log$-concave function, namely\footnote{$f\propto g$ means that
$f/g$ is constant.} $q(u,v)\propto \mathrm{e}^{-\alpha(u,v)}$, where $\alpha(u,v)$
is a non-negative twice continuously differentiable convex function on
$\Real^2$;

\item[(a3)]\label{a2} the observation density $p$ is a $C^2(\Real)$ $\log
$-concave function in the first argument: $ p(x,y)\propto \mathrm{e}^{-\gamma
(x,y)}$, where,
for each $y\in\Real$, the function $x\mapsto\gamma(x,y)$ is
non-negative, twice continuously differentiable and strongly convex on
$\Real$ with
$x_*(y):=\argmin_{x\in\Real}\gamma(x,y)\in(-\infty,\infty)$
and\vspace*{-2pt}
\[
\frac{\partial^2}{\partial x^2}\gamma(x,y)\ge\kappa>0 \qquad\forall x,y\in
\Real,
\]
with a constant $\kappa$;

\item[(a4)]\label{a3} for some constant $C$,
\[
-\mathop{\varlimsup}_{n\to\infty}\frac1 n \log L_n(X_{1\dvtx n},Y_{1\dvtx n}) \le C, \qquad\P\mbox{-a.s.;}
\]

\item[(a4)]\label{a4} there is a non-decreasing function $g\dvtx \Real_+\mapsto
\Real_+$, growing to $+\infty$ not faster than a polynomial, such that
for all $M>0$,
\[
\alpha(x,y)\le M \quad\implies\quad \bigg| \frac{\partial^2}{\partial x \partial y}\alpha(x,y) \bigg|\le g(M)\qquad \forall x,y\in\Real.
\]
\end{enumerate}
%r3.1
\begin{remark}
The $\log$-concavity assumptions (a1)--(a3) are quite restrictive.
For example, if $Y_n=h(X_n)+w_n$ with $w_n\sim N(0,1)$, then
\[
\frac{\partial^2} {\partial x^2}\gamma(x,y)=
\frac1 2 \frac{\partial^2} {\partial x^2} \bigl(y-h(x) \bigr)^2=
(h'(x) )^2- \bigl(y-h(x) \bigr)h''(x),
\]
which typically will not admit the uniform lower bound of (a3),
unless $h$ is linear, that is, $h''(x)\equiv0$.

If the assumption (a3) is satisfied, then it implies that
$\gamma_*(y):=\gamma(x_*(y),y )\in(-\infty,\infty)$ for all $y\in
\Real$ and, moreover,
%
%e11 ###
\begin{equation}
\label{kappa}
\gamma(x,y)- \gamma_*(y) \ge\tfrac1 2 \kappa(x-x_*)^2\qquad \forall x,y\in
\Real,
\end{equation}
which is essential to our approach.

Assuming that $-\log\mu(u)$, $\alpha(u,v)$ and $\gamma(x,y)$ are
non-negative is equivalent to assuming that they are lower-bounded by a
constant, that is, that the corresponding densities are bounded.

The assumption (a4) is typically satisfied if the state process
$X$ is positively recurrent (explicit recurrence tests can be found in
\cite{MT09}; see also \cite{GCJL00}).
Finally, (a5) is a technical assumption which is satisfied in
most models of practical
interest.
\end{remark}
%ex3.1
\begin{example}
All of the above assumptions are satisfied for the linear HMM
\begin{eqnarray*}
 X_n &=& a X_{n-1} + v_n,\qquad n\ge1, \\
 Y_n &=& b X_n + w_n,
\end{eqnarray*}
where $|a|<1$ and $b\ne0$ are constants and $v=(v_n)_{n\ge1}$ and
$w=(w_n)_{n\ge1}$ are independent
sequences of i.i.d.~random variables with
\[
X_0,v_n\sim f_v(x)\propto \mathrm{e}^{-|x|^{2+\delta}}
\quad\mbox{and}\quad
w_n\sim f_w(x)\propto \mathrm{e}^{- x^2 (1+c|x|^{\delta'})}
\]
for some $\delta\ge0$ and $\delta'\ge0$, $c\ge0$.
\end{example}
%t3.1
\begin{theorem}\label{theorem}
The limit in \eqref{hat} exists $\P$-a.s.
\end{theorem}

\begin{pf}
To keep the notation simple, we shall prove the convergence in \eqref
{hat} for $m=1$, that is, the limit $\lim_{n\to\infty}\hat X_1^n$
exists $\P$-a.s. As will be clear from the proof below, the same
arguments imply convergence of
$\lim_{n\to\infty}\hat X_i^n$ for any $i\le m$ and hence of \eqref{hat}
for any fixed $m\ge1$.

To check $\lim_{n\to\infty}\hat X_1^n$, $\P$-a.s., we shall show that
on a set of probability one, the series
\[
\hat X^n_1 = \hat X^1_1+\sum_{i=2}^n (\hat X^k_1-\hat X^{k-1}_1 )
\]
is convergent. The proof hinges on the system of inequalities \eqref
{mineqj} and \eqref{mineqn},
which stem from the $\log$-concavity properties assumed in (a1)--(a3). A pigeonhole principle
type of argument (Lemma~\ref{lemA1}) shows that a sequence satisfying
such inequalities must decay at least
polynomially backward in time, which, in turn, yields the desired conclusion.

To this end, introduce\footnote{For $k>\ell$, $\sum_{i=k}^\ell\cdots
=0$ is understood. }
%
%e12 ###
%
\begin{eqnarray}\label{hn}
h_n(x_{1\dvtx n}) &:=& -\log L_n(x_{1\dvtx n},Y_{1\dvtx n})\nonumber
\\[-8pt]\\[-8pt]
&=&-\log\mu(x_1) + \gamma(x_1,Y_1)+\sum_{m=2}^n \bigl(\alpha
(x_{m-1},x_m)+\gamma(x_m,Y_m) \bigr).\nonumber
\end{eqnarray}
By assumptions (a1)--(a3), $\lim_{R\to\infty}\inf_{\|
x_{1\dvtx n}\|= R}h_n(x_{1\dvtx n})\to\infty$ and, for any $n\ge1$, the function
%
%e13 ###
\begin{equation}\label{halpha}
x_{1\dvtx n}\mapsto h_n(x_{1\dvtx n})+\alpha(x_n,u)
\end{equation}
attains its global minimum at
\[
\tilde X^n_{1\dvtx n}(u) := \argmin_{x_{1\dvtx n}} \bigl(h_n(x_{1\dvtx n})+\alpha(x_n,u) \bigr),
\qquad u\in\Real.
\]
The Hessian matrix of the function defined in \eqref{halpha} is
positive definite uniformly over $x_{1\dvtx n}\in\Real^n$ and hence the
minimum is unique
and $\tilde X^n_{1\dvtx n}(u)$ is the solution of
\[
\operatorname{grad} \bigl(h_n(x_{1\dvtx n})+\alpha(x_n,u) \bigr)=0.
\]
The Jacobian matrix of the function on
the left-hand side of this equation with respect to the vector
$x_{1\dvtx n}$ coincides with the aforementioned Hessian matrix and
hence is invertible at any $u\in\Real$. Thus, by the implicit function
theorem, $u\mapsto\tilde X^n_{1\dvtx n}(u)$ is continuously differentiable
on $\Real$.

The usual dynamical programming argument yields the following chain rules:
%
%e14 ###
%
\begin{eqnarray}\label{chr1}
\tilde X_j^n(x) &=& \tilde X^m_j (\tilde X^n_{m+1}(x) ),\qquad x\in\Real, j<n, m=j,\ldots,n,\nonumber \\[-8pt]\\[-8pt]
\hat X^n_j &=& \tilde X^m_j (\hat X^n_{m+1} ).\nonumber
\end{eqnarray}
Hence, for $j< n$, and $j\le m < n$,
%
%e15 ###
%
\begin{eqnarray}\label{vottak}
\hat X^{n+1}_j -\hat X^n_j &=& \tilde X^{m}_j (\hat X^{n+1}_{m+1}
)-\tilde X^{m}_j (\hat X^n_{m+1} )\nonumber
\\[-8pt]\\[-8pt]
&=&
(\hat X^{n+1}_{m+1}-\hat X^n_{m+1} )\int_0^1 \frac\partial{\partial
s} \tilde X^{m}_j \bigl(s \hat X^{n+1}_{m+1} +(1-s)\hat X^n_{m+1} \bigr)\,\mathrm{d}s.\nonumber
\end{eqnarray}

The following lemma is the key to a bound on the integrand in \eqref{vottak}.
%l3.1
\begin{lemma}\label{lem34}
Assume \textup{(a1)--(a3)}. Then, for $j=1,\ldots,n-1$,
\begin{eqnarray}
\label{mineqj}
\bigg\|\frac{\partial}{\partial x}\tilde X^n_{1:j}(x) \bigg\|^2
\le
\frac2 {\kappa} \bigg|\Partial\alpha(\tilde X^n_{j}(x), \tilde
X^n_{j+1}(x) )
\frac\partial{\partial x}\tilde X^n_{j+1}(x)
\frac\partial{\partial x} \tilde X^n_{j}(x)
\bigg|%\succsim
\end{eqnarray}
and\vspace*{1pt}
\begin{eqnarray}
\label{mineqn}
\bigg\|\frac\partial{\partial x}\tilde X^n_{1\dvtx n}(x) \bigg\|^2 \le\frac2
{\kappa} \bigg| \Partial\alpha(\tilde X^n_{n}(x),x )
\frac\partial{\partial x} \tilde X^n_{n}(x) \bigg|,
\end{eqnarray}
where $\Partial\alpha(x,y):=\frac{\partial^2} {\partial x\,\partial
y}\alpha(x,y)$ and $\kappa$ is as in assumption (a3).
\end{lemma}

\begin{pf}
Recall that the function \eqref{halpha} is strongly convex and
the spectral norm of its Hessian is lower bounded by $\kappa$.
%Recall that the function \eqref{halpha} is convex with the Hessian,
%greater than $\kappa$ times the identity matrix, with respect to the
%positive definite ordering.
Hence, for any $1\le j<n$ and $u,v\in\Real$, by~\eqref{kappa},
\[
\frac{\kappa} 2 \|\tilde X^j_{1:j}(v)-\tilde X^j_{1:j}(u) \|^2\le
h_j (\tilde X^j_{1:j}(v) )+\alpha(\tilde X^j_{j}(v),u )-
h_j (\tilde X^j_{1:j}(u) )-\alpha(\tilde X^j_{j}(u),u )
\]
since,\vspace*{-3pt} by definition, the minimum of $h_j(x_{1:j})+\alpha(x_j,u)$ over
$x_{1:j}$ is attained at $\tilde X^j_{1:j}(u)$.
Further, by the definition of $\tilde X^j_{1:j}(v)$,
\[
h_j (\tilde X^j_{1:j}(v) )+\alpha(\tilde X^j_{j}(v), v )\le
h_j (\tilde X^j_{1:j}(u) )+\alpha(\tilde X^j_{j}(u), v ),
\]
which gives
%
%e16 ###
\begin{equation}
\label{yanki}
\frac{\kappa} 2 \|\tilde X^j_{1:j}(v)-\tilde X^j_{1:j}(u) \|^2\le
-\alpha(\tilde X^j_{j}(v), v )+\alpha(\tilde X^j_{j}(u),v )+\alpha
(\tilde X^j_{j}(v),u )
-\alpha(\tilde X^j_{j}(u),u ).
\end{equation}
Plugging $v: = \tilde X^n_{j+1}(x+h)$ and $u: = \tilde X^n_{j+1}(x)$
into this with $x\in\Real$ and using the chain rule \eqref{chr1}, we get
\begin{eqnarray*}
\frac{\kappa} 2 \|\tilde X^n_{1:j}(x+h)-\tilde X^n_{1:j}(x) \|^2&\le&
-\alpha\bigl(\tilde X^n_{j}(x+h), \tilde X^n_{j+1}(x+h) \bigr)
+\alpha\bigl(\tilde X^n_{j}(x),\tilde X^n_{j+1}(x+h) \bigr)
\\
&&{}+
\alpha\bigl(\tilde X^n_{j}(x+h),\tilde X^n_{j+1}(x) \bigr)
-\alpha(\tilde X^n_{j}(x),\tilde X^n_{j+1}(x) ).\vadjust{\goodbreak}
\end{eqnarray*}
Since all of the functions appearing in the latter inequality are twice
continuously differentiable, dividing by $h^2$ and taking $h\to0$ gives
the bound \eqref{mineqj}. Similarly, with $j:=n$, $v:=x+h$ and $u:=x$,
\eqref{yanki} yields \eqref{mineqn}.
\end{pf}

By assumption (a4),
\[
\Omega' := \Biggl\{ \mathop{\varlimsup}_{n\to\infty}\frac1 n \sum_{j=2}^n \bigl(\alpha
(X_{j-1},X_j)+\gamma(X_j,Y_j) \bigr)\le C \Biggr\}
\]
is an event of full probability and hence it is enough to verify the
claimed convergence for all $\omega\in\Omega'$.
Clearly, for an $\omega\in\Omega'$,
\[
-\log\mu(X_1) + \gamma(X_1,Y_1) +\sum_{j=2}^n \bigl(\alpha
(X_{j-1},X_j)+\gamma(X_j,Y_j) \bigr)\le2Cn\qquad \forall n\ge N(\omega)
\]
for an integer $N(\omega)<\infty$.
Then, $\hat X^n_{1\dvtx n}$, being a minimizer, a fortiori satisfies
%
%e17 ###
\begin{equation}
\label{ngeN}
-\log\mu(\hat X^n_1) + \gamma(\hat X^n_1,Y_1)+
\sum_{j=2}^n \bigl(\alpha(\hat X^n_{j-1},\hat X^n_j)+\gamma(\hat X^n_j,Y_j)
\bigr)\le2Cn\qquad \forall n\ge N.
\end{equation}
Hence, for a large fixed constant $M>4C$ and any $n\ge N$,
\[
\# \{j\dvtx  \alpha(\hat X^n_{j-1},\hat X^n_j)+\gamma(\hat X^n_j,Y_j) > M \}
\le \frac{2C n}{M}=: \rho n.
\]
Similarly,
\[
\# \{j\dvtx  \alpha(\hat X^{n+1}_{j-1},\hat X^{n+1}_j)+\gamma(\hat
X^{n+1}_j,Y_j) > M \}\le \frac{2C (n+1)}{M}= \rho(n+1).
\]
There is then an index $m\in[n-2\rho n,n]$ such that
\[
\alpha(\hat X^n_{m-1},\hat X^n_{m})+\gamma(\hat X^n_{m},Y_{m}) \le M
\quad\mbox{and}\quad
\alpha(\hat X^{n+1}_{m-1},\hat X^{n+1}_{m})+\gamma(\hat
X^{n+1}_{m},Y_{m}) \le M,
\]
and, by the assumption (a3),
\begin{eqnarray*}
|\hat X^{n+1}_{m}-\hat X^n_{m} |&\le&
\Big|\hat X^{n+1}_{m}-\argmin_{x\in\Real} \gamma(x,Y_{m}) \Big|
+ \Big|\hat X^{n}_{m}-\argmin_{x\in\Real} \gamma(x,Y_{m}) \Big|
\\
&\le&
\sqrt{\frac2 \kappa} \bigl(\gamma(\hat X^{n+1}_{m},Y_{m})- \gamma_* (Y_{m}) \bigr)^{1/2}+
\sqrt{\frac2 \kappa} \bigl(\gamma(\hat X^n_{m},Y_{m})- \gamma_* (Y_{m}) \bigr)^{1/2}
\\
&\le&\sqrt{\frac2 \kappa\gamma(\hat X^n_{m},Y_{m})}
+\sqrt{\frac2 \kappa\gamma(\hat X^{n+1}_{m},Y_{m})}\le\sqrt{\frac
{8M} \kappa}.
\end{eqnarray*}
Plugging this estimate into \eqref{vottak}, we get (for $j:=1$)
%
%e18 ###
\begin{equation}
\label{cbnd}
|\hat X^{n+1}_{1} -\hat X^n_{1} | \le
\sqrt{\frac{8M} \kappa}
\int_0^1 \bigg|\frac\partial{\partial s} \tilde X^{m-1}_{1} \bigl(s \hat
X^{n+1}_{m} +(1-s)\hat X^n_{m} \bigr) \bigg|\,\mathrm{d}s.
\end{equation}

Introduce
\begin{eqnarray*}
\check X^{m}_{m}(s)&:=&
s \hat X^{n+1}_{m} +(1-s)\hat X^n_{m},\\
\check X^{m}_j(s) &: =&
\tilde X^{m-1}_{j} (\check X^{m}_{m}(s)
),\qquad j =1,\ldots,m-1,
\end{eqnarray*}
and define
\begin{eqnarray*}
c_j(s) &:=& \frac2 \kappa |\Partial\alpha(\check X^{m}_j(s), \check
X^{m}_{j+1}(s) ) |,\qquad j < m,\\
b_j(s)&:=& \bigg|\frac\partial{\partial x} \tilde X^{m-1}_{j} (\check
X^{m}_{m}(x) )\Big|_{x:=s} \bigg|.
\end{eqnarray*}
Then, from \eqref{mineqj} and \eqref{mineqn} (the dependence on $s$ is
now omitted for brevity),
\begin{eqnarray}\label{bj}
\sum_{i=1}^j b_{i}^2&\le&
c_j
b_jb_{j+1}, \qquad j <m-1,\\
\label{bn}
\sum_{i=1}^{m-1} b_{i}^2 &\le& c_{m-1} b_{m-1}
\end{eqnarray}
and \eqref{cbnd} reads\vspace*{-1pt}
%
%e19 ###
\begin{equation}
\label{cbnd-reads}
|\hat X^{n+1}_{1} -\hat X^n_{1} | \le
\sqrt{\frac{8M} \kappa}
\int_0^1 b_1(s)\, \mathrm{d}s.
\end{equation}
%l3.2
\begin{lemma}
For any $s\in[0,1]$, $x>0$ and $g(\cdot)$ as in \textup{(a5)},\vspace*{-1pt}
%
%e20 ###
\begin{equation}
\label{yab}
\# \biggl\{j< m\dvtx  c_j(s) > \frac2 \kappa g(x) \biggr\}\le\frac{4C}{x(1-2\rho)} m.
\end{equation}
\end{lemma}
\begin{pf}
The function $u\mapsto\min_{x_{1\dvtx n}} (h_n(x_{1\dvtx n})+\alpha(x_n,u) )$ is
convex and hence\vspace*{-1pt}
\begin{eqnarray*}
\sum_{j=2}^{m} \alpha(\check X^{m}_{j-1},\check X_{j}^{m} ) &\le&
h_{m-1}(\check X^{m}_{1:m-1})+\alpha(\check X^{m}_{m-1},\check
X_{m}^{m} )
\\&=&
\min_{x_{1:m-1}} \bigl(h_{m-1}(x_{1:m-1})+\alpha\bigl(x_{m-1},s \hat
X^{n+1}_{m} +(1-s)\hat X^n_{m} \bigr) \bigr)
\\
&\le&
s \min_{x_{1:m-1}} \bigl(h_{m-1}(x_{1:m-1})+\alpha(x_{m-1}, \hat
X^{n+1}_{m} ) \bigr)
\\
&&{}+(1-s)\min_{x_{1:m-1}} \bigl(h_{m-1}(x_{1:m-1})+\alpha
(x_{m-1},\hat X^n_{m} ) \bigr)
\\
&=& s \bigl(h_{m-1}(\hat X^{n+1}_{1:m-1})+\alpha(\hat X^{n+1}_{m-1}, \hat
X^{n+1}_{m} ) \bigr)
\\
&&{}+(1-s) \bigl(h_{m-1}(\hat X^{n}_{1:m-1})+\alpha(\hat
X^{n}_{m-1},\hat X^n_{m} ) \bigr)
\\
&\le& 2C(n+1),
\end{eqnarray*}
where the latter inequality follows from \eqref{ngeN}. Hence,
\[
\# \{j\le m\dvtx  \alpha(\check X^{m}_{j-1},\check X_{j}^{m} )> x \}\le
\frac{2C(n+1)}{x}
\]
and, since $m\ge(1-2\rho)n$, \eqref{yab} follows from the assumption
(a5).
\end{pf}

Now, by Corollary \ref{corA} in the \hyperref[app-A]{Appendix}, applied to \eqref
{bj}--\eqref{bn} and \eqref{yab}, for any $\beta>1$,
there is a constant $C_\beta$ such that
%
%e21 ###
\begin{equation}
\label{bconv}
b_1 \le C_\beta m^{-\beta}\le C_\beta(1-2\rho)^{-\beta}n^{-\beta}
\end{equation}
for all sufficiently large $n$ and, thus, by \eqref{cbnd-reads}, the
sequence $ |\hat X^{n+1}_{1} -\hat X^n_{1} |$, $n\ge1,$ is
summable, which verifies the existence of the limit \eqref{hat}.
\end{pf}

\begin{remark}\label{rem-G}
When the hidden state process is a Gaussian autoregression, that is,
when $\alpha(x,y)=\frac1 2(y-bx)^2$ with a constant $b\ne0$, $
|\Partial\alpha(x,y) |\equiv b$
and Lemma \ref{lemA1}(1) implies the exponential bound in \eqref
{bconv}, confirming the results deducible from Kalman linear filtering theory.
\end{remark}

%s4 ###
\section{Concluding remarks}\label{sec-4}
As indicated by the examples of Section \ref{sec-2} and the partial
results of Theorem \ref{theorem},
the convergence in \eqref{hat} appears to be a non-trivial issue.
Analogous problems have been discussed in
the engineering literature.
In fact, the MAP path estimation can be viewed as an optimal control
problem, in which
one is required to minimize the cost functional $h_n(x_{1\dvtx n})$ defined
in \eqref{hn}, where the term $\alpha(x_{m-1},x_m)$ is interpreted as
the cost
incurred by the control effort (needed to move from $x_{m-1}$ to $x_m$)
and $\gamma(x_m,Y_m)$ is
the cost paid for the deviation of the state from $Y_m$.
This setting appears in \cite{Bel}, Chapter I,
Section 1.7, as the ``smoothing'' problem and,
in the control literature, is often referred to as the \textit
{tracking} problem.
From the control theory perspective, the existence of the limit in
\eqref{hat} means that the optimal control and the corresponding
optimal trajectory cease to depend on the future values of the
exogenous signal $Y$.

Among other related questions, the convergence \eqref{hat} of the
optimal trajectory is part of the
``asymptotic control theory'' program initiated by R.~Kalman,
R.~Bellman and R.~Bucy, at the dawn of
modern control theory.
In the linear state/quadratic cost (LQ) setting of R.~Kalman, the
control problem admits an elegant closed-form
solution for each fixed horizon $n$ and the study of the limit \eqref
{hat} reduces to the stability analysis of
the associated Riccati equation (a~comprehensive treatment of the LQ
problem can be found in, e.g., \cite{KS72}).

To the best of our knowledge, asymptotic analysis beyond the LQ case
has been carried out only for a limited number of nonlinear
models. Bellman and Bucy \cite{BB64} found a remarkable explicit solution to a quite
general scalar continuous-time control problem, amenable to asymptotic
analysis. A vector control problem with linear state dynamics and
convex costs was studied in~\cite{B66}.

While much progress has been made in the optimal control theory on the
\textit{infinite horizon} (see, e.g., \cite{CHL91,Za06}),
we were not able to track any results directly applicable to the
question under consideration.

Another possible connection, remaining elusive at the moment, is to the
stability theory of nonlinear filtering equations,
developed during the last decade (see, e.g., the survey~\cite{CLV09}).

\begin{appendix}
%s5 ###
\section*{Appendix: A supporting lemma}\label{app-A}

\begin{lemma}\label{lemA1}
Consider the system of inequalities
%
%e22 ###
\begin{eqnarray}\label{ineq}
\sum_{i=1}^j b_i^2&\le& b_jb_{j+1} c_j,\qquad j=1,\ldots,n-1,\nonumber\\[-8pt]\\[-8pt]
\sum_{i=1}^n b_i^2&\le& b_n c_n,\nonumber
\end{eqnarray}
where $b_i$ and $c_i$, $i=1,\ldots,n,$ are non-negative real numbers,
and let $\theta$ and $\theta'$ be
arbitrary positive constants:
\begin{enumerate}[(1)]
\item[(1)] If $c_i \le\theta$, $i=1,\ldots,n$, then
%
%e23 ###
\begin{equation}
\label{b1s}
b_1\le
\sqrt{\theta \mathrm{e}} \exp\biggl(-\frac{n}{2\mathrm{e}(\theta^2\vee\theta)} \biggr)\qquad \mbox{for }
n\ge\theta^2 \mathrm{e}.
\end{equation}

\item[(2)] If, for a non-decreasing non-negative function $g\dvtx \Real_+\mapsto
\Real_+$,
%
%e24 ###
\begin{equation}
\label{MI}
\#\{i\le n\dvtx  c_i \ge g(x)\}\le\frac{\theta n}{x}\qquad \forall x>0,
\end{equation}
and $c_n \le\theta'$, then, for any $p\in(0,1)$ and $\ell> \theta$,
%
%e25 ###
\begin{equation}
\label{b1}
b_1\le\sqrt{g(\ell)} n^{-p\ell/(4\theta)}\qquad \mbox{for } n > \biggl(\frac{\ell
(\theta'^2\vee g(\ell))}{\theta} \biggr)^{1/(1-p)}.
\end{equation}

\item[(3)]\label{claim3} If only \textup{\eqref{MI}} holds, then, for any $p\in
(0,1)$ and $\ell> \theta$,
%
%e26 ###
\begin{equation}
\label{b1cor}
b_1\le g(2\theta n)\sqrt{g(\ell)} n^{-p\ell/(4\theta)}\qquad \mbox{for } n >
\biggl(\frac{\ell(1\vee g(\ell))}{\theta} \biggr)^{1/(1-p)}.
\end{equation}
\end{enumerate}
\end{lemma}

\begin{pf}
(1) The second inequality in \eqref{ineq} and $c_n\le\theta$ together
imply that $b_n^2\le b_n\theta$ and, in turn, that $b_1^2+ \cdots+
b_n^2\le \theta^2$.
Fix a constant $\eta\in(0,1)$ and let $m_1:=\lfloor\theta^2/\eta
\rfloor$. Then, at most half of the $b_i$'s with $i\in[n-2m_1,n]$ are
greater than
$\sqrt{\eta}$ and hence there is an index $k_1 \in[n-2m_1,n]$ such
that $b_{k_1}\le\sqrt{\eta}$ and $b_{k_1+1}\le\sqrt{\eta}$.
The inequality corresponding to $j:=k_1$ in \eqref{ineq} then gives the
bound $b_1^2+\cdots+b_{k_1}^2\le b_{k_1} b_{k_1+1}c_{k_1}\le\eta\theta$.

Similarly, let $m_2:=\lfloor\theta/\eta\rfloor$. There is then an
index $k_2 \in[k_1-2m_2:k_1]$ such that $b_{k_2}\le\eta$ and
$b_{k_2+1}\le\eta$ and, again applying \eqref{ineq}, $b_1^2+\cdots+
b_{k_2}^2\le b_{k_2}b_{k_2+1}c_{k_2}\le\eta^2\theta$. This argument
can be iterated at least
\[
\bigg\lfloor
\frac{n}{2(m_1 \vee m_2)}
\bigg\rfloor=
\bigg\lfloor
\frac{\eta n}{2(\theta^2 \vee\theta)}
\bigg\rfloor
\]
times and thus
\[
b_1^2 \le\theta\eta^{ \lfloor
\sklfrac1 2\eta n/(\theta^2 \vee\theta)
\rfloor}\le
\frac{\theta}{\eta} \bigl((\eta^\eta)^{\sklfrac1 2 /(\theta^2\vee\theta)} \bigr)^{n}.
\]
The best rate is obtained at $\eta:=\mathrm{e}^{-1}$, which yields the bound
\eqref{b1s}.

(2) For a fixed $\ell\ge\theta$, by \eqref{MI},
%
%e27 ###
\begin{equation}
\label{sharp}
\#\{i\le n\dvtx  c_i \ge g(\ell)\}\le\frac{\theta n}{\ell}:=r n
\end{equation}
and thus at least half of the $c_i$'s with $i\in[n-2rn,n]$ do not
exceed $g(\ell)$.
Fix a constant $p\in(0,1)$ and let $\eta:=n^{-p/2}$. Suppose that for
all $i\in[n-2rn,n]$ such that $c_i\le g(\ell)$, either $b_i\ge\eta$
or $b_{i+1}\ge\eta$, or both. Then,
\[
\#\{i\in[n-2rn\dvtx n]\dvtx  b_i\ge\eta\}\ge rn.
\]
But, on the other hand, by the second inequality in \eqref{ineq} and as
$c_n\le\theta'$, $b_n^2\le b_n \theta'$ and
$b_1^2+\cdots+b_n^2\le\theta'^2$, we have
\[
\#\{i\in[n-2rn\dvtx n]\dvtx  b_i\ge\eta\}\le\frac{\theta'^2}{\eta^2}=\theta'^2 n^{p}.
\]
This contradicts the previous estimate if $n$ is large enough, namely,
if $n> (\ell\theta'^2/\theta)^{1/(1-p)}$. Thus,
for such $n$, there is an index $m_1\in[n-2rn\dvtx n]$ such that
$b_{m_1}\le\eta$, $b_{m_1+1}\le\eta$ and $c_{m_1}\le g(\ell)$.

Now, by the inequality in \eqref{ineq} corresponding to $j:=m_1$,
%
%e28 ###
\begin{equation}
\label{etaell}
b_1^2+\cdots+b_{m_1}^2\le b_{m_1}b_{m_1+1}c_{m_1}\le\eta^2g(\ell)
\end{equation}
for which the above consideration can be repeated. Namely,
by \eqref{sharp}, there are at least $rn$ indices $i\in[m_1-2rn,m_1]$
for which $c_i\le g(\ell)$.
Suppose that, for all of them, either $b_i \ge\eta^2$ or $b_{i+1}\ge
\eta^2$, or both. Then
\[
\#\{i\in[m_1-2rn\dvtx m_1]\dvtx  b_i\ge\eta^2\}\ge rn,
\]
while \eqref{etaell} implies
that
\[
\#\{i\in[m_1-2rn\dvtx m_1]\dvtx  b_i\ge\eta^2\}\le\frac{\eta^2g(\ell)}{\eta
^4}=n^p g(\ell),
\]
which is a contradiction for $n$ large enough, that is, for $n> (\ell
g(\ell)/\theta)^{1/(1-p)}$. Hence, there is an $m_2\in[m_1-2rn\dvtx m_1]$
such that
$b_{m_2}\le\eta^2$, $b_{m_2+1}\le\eta^2$ and $c_{m_2}\le g(\ell),$
and thus, by \eqref{ineq}, we have
\[
b_1^2+\cdots+b_{m_2}^2 \le\eta^4g(\ell).
\]
This argument can be iterated at least $\lfloor1/(2r)\rfloor$ times,
which yields the bound
\[
b_1^2 \le g(\ell) (\eta^{1/2r} )^2 = g(\ell) n^{-p\ell/2\theta}.
\]

(3) Note that $b'_i := b_i/g(2\theta n)$, $i=1,\ldots,n$, satisfy the
inequalities \eqref{ineq} with the $c_i$'s replaced by $c'_i := c_i$,
$i=1,\ldots,n-1$, and
$c'_n:=c_n/g(2\theta n)$. By \eqref{MI},\vspace{-5pt}
\[
\# \{i\le n\dvtx  c_i \ge g(2\theta n) \} \le\frac{\theta n}{2\theta n}=1/2,
\]
that is, all $c_i$'s are less than $g(2\theta n)$ and, in particular,
$c_n\le g(2\theta n)$, that is, $c'_n\le1$. Moreover, assuming that
$g(2\theta n)\ge1$,
we have
\[
\# \{i\le n\dvtx  c'_i \ge g(x) \}\le\# \{i\le n: c_i \ge g(x) \} \le\frac
{\theta n}{x}\qquad \forall x>0.
\]
Hence, by \eqref{b1}, we have\vspace{-5pt}
\[
b'_1\le\sqrt{g(\ell)} n^{-p\ell/(2\theta)}\qquad \mbox{for } n > \biggl(\frac{\ell
(1\vee g(\ell))}{\theta} \biggr)^{1/(1-p)},\vspace{-5pt}
\]
which, in turn, gives \eqref{b1cor}.
\end{pf}
%cA.1
\begin{coro}\label{corA}
Under the assumption \textup{\eqref{MI}} with $g(\cdot)$ growing to $+\infty$
not faster than a polynomial,
for any $\beta>1$, there is a constant $C_\beta$, such that\vspace{-5pt}
\[
b_1 \le C_\beta n^{-\beta}\vspace{-5pt}
\]
for all sufficiently large $n$.
\end{coro}
\begin{pf}
This follows from (3) %\eqref{claim3}
of Lemma \ref{lemA1}.
\end{pf}
\end{appendix}

\section*{Acknowledgements} The authors are grateful to M.~Margaliot
and R.~Van Handel for the enlightening discussion
on the issues raised in this paper and appreciate the thorough
proofreading by the referees. Research of P.~Chigansky was supported by
an ISF Grant 314/09. Research of Yaacov Ritov was supported by an ISF grant.

\printhistory

\end{document}